\documentclass[11pt, reqno]{amsart}

 \usepackage{amsmath}
 \usepackage{amssymb}
\usepackage{bm}

\newcommand{\mysection}[1]{\section{#1}
      \setcounter{equation}{0}}

\newtheorem{theorem}{Theorem}[section]
\newtheorem{lemma}[theorem]{Lemma}

\theoremstyle{definition}
\newtheorem{assumption}{Assumption}[section]

\theoremstyle{remark}
\newtheorem{remark}{Remark}[section]

\newcommand{\tr}{\text{\rm tr}\,}

\newcommand\bbeta{\text{\raise-.2ex\hbox{$\bm{\beta}$}}}
\newcommand\dist{{\rm dist}\,}

\newcommand\bR{\mathbb{R}}

\newcommand\bS{\mathbb{S}}

\newcommand\bZ{\mathbb{Z}}

\newcommand\cL{\mathcal{L}}

\newcommand\supinf{\operatornamewithlimits{sup\,\,\,inf}}

\begin{document}

\title[Isaacs equations]
{To the theory of viscosity solutions
for uniformly elliptic  Isaacs equations}

\author{N.V. Krylov}
\thanks{The  author was partially supported by
 NSF Grant DMS-1160569}
\email{krylov@math.umn.edu}
\address{127 Vincent Hall, University of Minnesota,
 Minneapolis, MN, 55455}



\begin{abstract}
We show how a theorem about  the
solvability in $C^{1,1}$ of special Isaacs equations
can be used to obtain existence and uniqueness
of viscosity solutions of general
uniformly nondegenerate Isaacs equations.
We apply it also to establish the $C^{1+\chi}$
regularity of viscosity solutions and
show that finite-difference approximations
have an algebraic rate of convergence.
The main coefficients of the Isaacs equations
are supposed to be in $C^{\gamma}$ with $\gamma$
slightly less than $1/2$.
\end{abstract}

\maketitle

\mysection{Introduction}
                                           \label{section 4.3.1}

The goal of this article is to present
a purely PDE exposition of some major
results in the theory of viscosity
solutions for uniformly
nondegenerate Isaacs equations.

Let $\bR^{d}=\{x=(x^{1},...,x^{d})\}$
be a $d$-dimensional Euclidean space.
Assume that we are given separable metric spaces
  $A$ and $B$,   and let,
for  $(\alpha,\beta,x)\in A\times B\times\bR^{d}$, 
  the following  
  functions  be given: 

(i) $d\times d $
matrix-valued $a^{\alpha\beta}( x)$,

(ii)
$\bR^{d}$-valued $b^{\alpha\beta}( x)$, and

(iii)
real-valued  functions 
$c^{\alpha\beta}( x) \geq0$,   
  $f^{\alpha\beta}( x) $, and  
$g(x)$. 

Let $\bS$ be the
set of symmetric $d\times d$ matrices, and for
$(u_{ij})\in\bS$, $(u_{i})\in\bR^{d}$, and $u\in\bR$
introduce 
$$
F(u_{ij},u_{i},u,x)=
\supinf_{\alpha\in A\,\,\beta\in B}
[a^{\alpha\beta}_{ij}( x)u_{ij}  +
b ^{\alpha\beta}_{i }( x)u_{i} -c^{\alpha\beta} ( x)u],
$$
where and everywhere below the summation convention is enforced
and the summations are done inside the brackets.

For a sufficiently smooth function $u=u(x)$ also introduce
$$
L^{\alpha\beta} u(x)=a^{\alpha\beta}_{ij}( x)D_{ij}u(x)+
b ^{\alpha\beta}_{i }( x)D_{i}u(x)-c^{\alpha\beta} ( x)u(x),
$$
where, naturally, $D_{i}=\partial/\partial x^{i}$, $D_{ij}=D_{i}D_{j}$.
Denote
\begin{equation}
                                                    \label{1.16.1}
F[u](x)=F(D_{ij}u(x),D_{i}u(x),u(x),x)
=\supinf_{\alpha\in A\,\,\beta\in B}
[L^{\alpha\beta} u(x)+f^{\alpha\beta} (x)].
\end{equation}

Also fix a sufficiently regular domain $G\subset \bR^{d}$.
Under appropriate assumptions which we list in
Section \ref{section 2.26.3} and which include the boundedness
and continuity with respect to $x$ of the data and uniform
nondegeneracy of $a^{\alpha\beta}( x)$ the Isaacs equation
\begin{equation}
                                                     \label{4.1.4}
F[u]=0
\end{equation}
in $G$ with boundary condition $u=g$ on $\partial G$
has a   viscosity solution $w\in C(\bar{G})$.
Recall (see \cite{CIL})
that this means that for any smooth $\phi(x)$
and any point $x_{0}\in G$ at which $\phi-w$ attains

(i) 
a local maximum which is zero we have $F[\phi](x_{0})\leq0$,

(ii) a local minimum
 which is zero we have $F[\phi](x_{0})\geq0$.

We are going to discuss the existence, uniqueness, regularity
properties of $w$, and the rate of convergence of finite-difference
approximations to $w$ and, therefore, we give 
 a brief account of basic facts known
for the Isaacs equations. We only discuss these equations
although in the references below more general equations
are considered and more details can be found.
For brevity, when we mention that, say, $a$ is uniformly
continuous in $x$, we mean uniformity with respect to $\alpha,
\beta,x$. The Lipschitz or H\"older continuity also presume
that the 
corresponding constants are independent of $\alpha,\beta$.
 
In 1989 Ishii \cite{Is_89} proved the existence
of viscosity solutions
for possibly degenerate equations
with Lipschitz continuous $a$ and $b$
and uniformly
continuous $c$ and $f$ with respect to $x$ under the condition
that for a constant $\delta>0$
\begin{equation}
                                                       \label{4.3.1}
c^{\alpha\beta}\geq\delta.
\end{equation}

In the same year 1989 Trudinger \cite{Tr_89} 
(see Corollary 3.4 there, also see \cite{Tr_88}) proved
the existence (for uniformly nondegenerate case)
dropping \eqref{4.3.1} and reducing the Lipschitz continuity
of $a$
to the $1/2-\varepsilon$-H\"older continuity and assuming that
$b,c,f$ are uniformly continuous with respect to $x$.
In \cite{Tr_89} the methods of the theory of viscosity 
solutions are combined with the classical PDE methods.

  Crandall,   Ishii, and Lions in their user's guide \cite{CIL}
in 1992 based their existence results
on the comparison principle in the same way as
it was done in \cite{Is_89} and one can extract from
\cite{CIL} an existence result
(unfortunately not stated explicitly) in the uniformly
nondegenerate  case 
under, basically, the same assumptions as in \cite{Is_89}
but dropping \eqref{4.3.1}.

Jensen and  \'Swi{\c e}ch  \cite{JS_05}  in 2005
proved that even if $a,b,c,f$ are just measurable,
there exists a continuous $L_{p}$-viscosity solutions,
which in case that
$a, b,c,f$ are uniformly continuous with respect to $x$
is automatically just a classical viscosity
 solution. To the best of the author's knowledge, 
this provides the most general
conditions to date for existence of classical viscosity
solutions (their solution in case of continuous $a$
 is actually in 
  $C^{1+\chi}$ which is proved in \'Swi{\c e}ch \cite{Sw97} (1997)
 and
can be obtained from an elliptic
counterpart of \cite{Kr_14}).
For further information
concerning $L_{p}$-viscosity solutions
we refer the reader to
\cite{CCKS_96} and \cite{CKS00}

Next issue is uniqueness. Here the starting
assumption is that the coefficients and $f$
are uniformly continuous in $x$.
Jensen \cite{Je_88} in 1988 proved a comparison principle 
(and hence uniqueness)
for {\em Lipschitz\/} continuous viscosity
solutions to the fully nonlinear second 
order elliptic PDE not explicitly depending on
$x$ and not necessarily connected with the Isaacs equations.
 Ishii in \cite{Is_89}
showed among other things that the Lipschitz continuity
of solutions is not necessary and one can treat
equations with $F$ explicitly depending on
$x$. The equations in \cite{Is_89} could degenerate
and therefore condition \eqref{4.3.1} is imposed.
Ishii and Lions \cite{IL_90} (1990) prove uniqueness
when $a$ is $1/2+\varepsilon\,$-H\"older
continuous.
Even though the equations in \cite{IL_90}
 are uniformly nondegenerate
condition \eqref{4.3.1} is imposed.

Trudinger \cite{Tr_89} (1989)
shows that if $a$ is H\"older continuous,
then any continuous viscosity solution
is Lipschitz continuous
and then uniqueness (but not existence) follows from
his result in \cite{Tr_88} (1988)
if $a$ is $1/2-\varepsilon\,$-H\"older continuous
(see Remark 3.1 in \cite{Tr_88}).
To the best of the author's knowledge, 
this is the most general result concerning
the uniqueness of classical viscosity solutions.

 Crandall,   Ishii, and Lions   \cite{CIL}
stated and proved the comparison
result  under, basically, the same assumptions as in \cite{Is_89},
and then made a few comments about the possibility to drop
\eqref{4.3.1} and prove the comparison result
for subsolutions and supersolutions in class $C^{1+\chi}$.
These comments will be instrumental in our exposition.

Jensen and  \'Swi{\c e}ch  \cite{JS_05}  proved uniqueness
of continuous $L_{p}$-viscosity solutions,
allowing
$1/2-\varepsilon\,$-H\"older continuous $a$ and
only measurable $b,c,f$. Again, if $b,c,f$
are uniformly continuous in $x$, this yields
the uniqueness of classical viscosity solutions,
previously obtained in \cite{Tr_89}.

Next issue is about the regularity of viscosity
solutions: given a continuous viscosity solution $w$,
is it true that $w\in C^{\chi}$ or $C^{1+\chi}$?

Caffarelli \cite{Ca89}
(1989) and Trudinger \cite{Tr_88}, \cite{Tr_89}
were the first authors who proved $C^{1+\chi}$
regularity for fully nonlinear elliptic equations of type 
\begin{equation}
                                              \label{4.4.1}
F[u]=f
\end{equation}
without convexity assumptions on $F$. The assumptions in these 
papers are different. We keep concentrating
on the Isaacs equations and compare the assumptions only
in that case.
In \cite{Ca89} the function $F(u_{ij},u_{i},u,x)$
 is independent of $u_{i},u$  
and, for each $u_{ij}$, is uniformly sufficiently close to a 
function which is
 continuous
with respect to $x$. In \cite{Tr_88} and 
\cite{Tr_89} the function $F$ depends on all arguments but is
H\"older continuous in $x$. Next step in what concerns
$C^{1+\chi}$-estimates   was done by   
\'Swi{\c e}ch \cite{Sw97} (1997),
 who considered general $F$ and imposed
the same condition as in \cite{Ca89} on the $x$-dependence,
which is  much weaker than in \cite{Tr_88} and \cite{Tr_89} (but
also imposed the Lipschitz condition on the dependence of $F$ on
$u_{i},u$ for $F$  more general than coming 
from the Isaacs equations,
whereas
 in \cite{Tr_88} and \cite{Tr_89} only the 
continuity with respect to
$u_{i}$ is assumed). It is worth emphasizing
that these results are about {\em any\/} continuous viscosity
solution. The existence of {\em a\/} viscosity
solution of class $C^{1+\chi}$ 
follows from \cite{JS_05}, \cite{Sw97} and also
can be obtained from an elliptic
counterpart of \cite{Kr_14}.

Finally, a few words about the rate of
convergence of numerical approximation
of solutions of the Isaacs equations.
Caffarelli and  Souganidis in \cite{CS} proved
that there is an algebraic rate
of convergence of solutions of finite-difference
schemes to the Lipschitz continuous viscosity solution
of the fully nonlinear elliptic equation \eqref{4.4.1}
with $F$ not necessarily connected with the Isaacs equations
in a regular domain with Dirichlet boundary data.
They assumed that $F$ depends only on $u_{ij}$.
This was the first result available for
 fully nonlinear elliptic equations without
convexity assumptions on $F$. Naturally,
one would want to extend the result to the  $F$'s
depending also on $Du$, $u$, and $x$.
Turanova \cite{Tu_13_1} extended the results of 
\cite{CS} to $F$'s explicitly depending on $x$,
but still independent of $u_{i},u$
and, in case of Isaacs equations,
with Lipschitz continuous $a$ and $f$.
This result was generalized in \cite{Kr_14_2}
for equations with Lipschitz continuous $a$ and $b$
and H\"older continuous $c$ and $f$.

In this paper we further generalize the result of
\cite{Kr_14_2} in the case of $1/2-\varepsilon\,$-H\"older
continuous $a$ and H\"older
continuous $b,c,f$.
This is done by using  special approximations
of the original equation introduced in \cite{Kr_12}.
On the way we also prove the uniqueness
of viscosity solutions for $1/2-\varepsilon\,$-H\"older
continuous $a$ and continuous $c,b,f$
(Trudinger's result of 1989).
In that case we
 also prove the existence of viscosity solution
in class $C^{1+\chi}$. 
Actually, as it has been mentioned
in the Abstract and in the beginning of the Introduction,
the main goal of this article 
is to show how a theorem  from \cite{Kr_12} about  
the solvability in $C^{1,1}$
of special Isaacs equations
can be used to obtain existence and uniqueness
of viscosity solutions, their $C^{1+\chi}$
regularity, and establish a
 rate of convergence of numerical approximation for general
uniformly nondegenerate Isaacs equations.

The methods we use
are different from the methods of
Trudinger and the methods of the theory
of viscosity solutions. However, apart from 
the results of \cite{Kr_12}
and other PDE tools we also use an argument
from Section V.~A of \cite{CIL} explaining how to
prove the comparison principle for
$C^{1+\chi}$ subsolutions and supersolutions.
We use a quantitative version of this argument.

The article is organized as follows.
In Section \ref{section 2.26.3} we present
our main results and prove all of them
apart from Theorem \ref{theorem 4.1.2}
 and assertion (ii) of Theorem \ref{theorem 4.1.3},
which are proved in Sections \ref{section 4.3.3}
and \ref{section 4.3.4}, respectively, after
a rather long Section
\ref{section 4.3.2} containing a comparison
theorem for smooth functions.

Our equation are considered in $C^{2}$ domains
with $C^{1,1}$ boundary data. These
restrictions can be
considerably relaxed and we leave doing that to the interested
reader.
 
\mysection{Main result }
                                           \label{section 2.26.3}

Fix some constants $\delta\in(0,1)$ and $K_{0}\in[0,\infty)$.
 Set
$$
\bS_{\delta}=\{a\in\bS:\delta|\xi|^{2}\leq a_{ij}\xi^{i}\xi^{j}
\leq\delta^{-1}|\xi|^{2},\quad \forall\, \xi\in\bR^{d}\}. 
$$

In the following assumption the small parameter  $\chi\in(0,1)$,
 which depends only on
$\delta$ and $d$, is a constant  to be specified in
Theorem \ref{theorem 3.31.1} and
\begin{equation}
                                                      \label{4.14.2}
\gamma=\frac{4-3\chi}{8-4\chi}\quad(<1/2).
\end{equation}
 
\begin{assumption}
                                    \label{assumption 1.9.1}
(i) The   functions 
$a^{\alpha\beta}( x )$,
$b^{\alpha\beta}( x )$, $c^{\alpha\beta}( x )$,
and $f^{\alpha\beta}( x )$
are continuous with respect to
$\beta\in B$ for each $(\alpha, x)$ and continuous with respect
to $\alpha\in A$ uniformly with respect to $\beta\in B$
for each $x$, and  
$$
\|g\|_{ C^{1,1}(\bR^{d})}\leq K_{0},
$$ 

(ii) 
For any $x \in\bR^{d}$ 
 and $(\alpha,\beta )\in A\times B $
$$
  \|a^{\alpha\beta}( x )\|,|b^{\alpha\beta}( x )|,
|c^{\alpha\beta}( x )|,|f^{\alpha\beta}( x )|
\leq K_{0},
$$
where     for a matrix $\sigma$ we denote $\|\sigma\|^{2}
=\tr\sigma\sigma^{*}$,

(iii) For any $(\alpha,\beta )\in A\times B $ and $x,y\in\bR^{d}$
we have
$$
\|a^{\alpha\beta}( x )-a^{\alpha\beta}( y )\|
\leq K_{0}|x-y|^{\gamma},
$$
$$
|u^{\alpha\beta}( x )-u^{\alpha\beta}( y )| 
\leq K_{0}\omega(|x-y|),
$$
where   $u=b,c,f$, and $\omega$ is a fixed continuous
increasing 
function on $[0,\infty)$ such that $\omega(0)=0$,

(iv) For all values of arguments $a^{\alpha\beta}\in\bS_{\delta}$.
 
\end{assumption}

\begin{remark}
                                \label{remark 4.5.1}
It is convenient to assume that
$$
|x-y|^{\chi}\leq\omega(|x-y|),
$$
whenever $|x-y|\leq1$.
 Clearly this assumption
does not restrict generality.
\end{remark}

We will be dealing with  equation \eqref{4.1.4}
in a fixed bounded domain $G\in C^{2}$ with boundary condition $u=g$
on $\partial G$. 
In \cite{Kr_12}
 a convex positive homogeneous of degree
one Lipschitz continuous function $P(u_{ij},u_{i},u)$ is constructed
on $\bS\times\bR^{d}\times\bR$
such that at all points of differentiability
of $P$ with respect to $(u_{ij})$ we have 
$$(P_{u_{ij}})\in\bS_{\hat \delta},\quad
|(P_{u_{i}})|\leq K_{1},\quad 0>P_{u}\geq -K_{1},
$$
where $\hat{\delta}$ is a constant in $(0,\delta)$
depending only on
$d $  and $\delta$ and $K_{1}\geq K_{0}$ 
depends only on
$d $, $K_{0}$,  and $\delta$.  This function
is constructed once only $d$, $K_{0}$, and $\delta$ are given
and possesses some additional properties to be mentioned and used
below. By $P[u](x)$ we denote $P(D_{ij}u(x),D_{i}u(x),u(x))$.

Here is a consequence of Theorems 1.1 and 1.3 of \cite{Kr_12}.
\begin{theorem}
                                           \label{theorem 3.30.1}
For any $K\geq0$ each of the equations
\begin{equation}
                                                 \label{3.30.6}
\max(F[u],P[u]-K)=0,
\end{equation}
\begin{equation}
                                                 \label{3.30.7}
 \min(F[v],-P[-v]+K)=0
\end{equation}
in $G$ with boundary condition $u=v=g$ on $\partial G$
has a unique solution in the class $C^{1,1}_{loc}(G)
\cap C(\bar G)$. 

\end{theorem}

By $u_{K}$ and $v_{K}$ we denote
 the solutions of \eqref{3.30.6} and \eqref{3.30.7},
respectively. These are the central objects of our
investigation. Here is a simple property they
possess.

\begin{lemma}
                                       \label{lemma 4.1.1}
There exists a constant $N$, depending only on
$d,\delta$,  $K_{0}$, and $G$,  such that in $G$
$$
|u_{K}-g|+|v_{K}-g|\leq N\rho,\quad
|u_{K} |+|v_{K} |\leq N,
$$
where $\rho(x)=\dist(x,G^{c})$.
\end{lemma}

This result for $u_{K}$
follows from the fact that 
  $|\max(F[0],-K)|\leq|F[0]|$, $g\in C^{1,1}$, $G\in C^{2}$, and
$u_{K}$ satisfies a linear equation
$$
a_{ij}D_{ij}u_{K}+b_{i}D_{i}u_{K}-cu_{K}+f=0,
$$
where $(a_{ij})\in\bS_{\hat \delta}$, $|(b_{i})|\leq K_{1}$,
$K_{1}\geq c\geq 0$, and $|f|\leq K_{0}$.
The case of $v_{K}$ is quite similar.

To characterize some smoothness properties of $u_{K}$
and $v_{K}$
introduce $C^{1+\chi}(G)$ as the space 
of continuously differentiable functions in $G$
with finite norm given by
$$
\|u\|_{C^{1+\chi}(G)}=\sup_{G}|u|
+\sup_{x,y\in G}\frac{|u( x)-u( y)|}{|x-y|}
+[u]_{C^{1+\chi}(G)},
$$
where
$$
[u]_{C^{1+\chi}(G)}=
\sup_{x,y\in G} \frac{|Du( x)-Du( y)|}{|x-y|^{\chi}}.
$$

For $\varepsilon>0$ introduce
\begin{equation}
                                                 \label{4.6.1}
  G_{\varepsilon}=\{x\in G:
\dist (x,\partial G)>\varepsilon\}.
\end{equation}

\begin{theorem}
                                            \label{theorem 3.31.1}
There exists a constant $\chi\in(0,1)$,
depending only on
$\delta$ and $d$, and there exists a constant 
$N$, depending only on $K_{0}$,
$\delta$, $d$, and $G$, such that for any $\varepsilon\in(0,1]$
(such that $G_{\varepsilon}\ne\emptyset$)
\begin{equation}
                                                 \label{3.30.8}
 \|u_{K},v_{K}\|_{C^{1+\chi}(G_{\varepsilon})}\leq
N\varepsilon^{-1-\chi}.
\end{equation}

\end{theorem}

Proof. Denote 
$$
F_{K}[u]=\max(F[u],P[u]-K) , \quad P_{0}(u_{ij})
=P(u_{ij},0,0),
$$
$$
 P_{0}[u](x)=P_{0}(D_{ij}u(x)),
$$
take the constant
$N$ from Lemma \ref{lemma 4.1.1}, and 
for $K\geq K_{1}N$ consider the {\em parabolic\/} equation
\begin{equation}
                                                 \label{4.1.2}
 \partial_{t}u+\max(F_{K}[u],P_{0}[u]-K_{1}N-K)=0
\end{equation}
in $(0,1)\times G_{\varepsilon/2}$ with boundary condition
$u=u_{K}$ on the parabolic boundary of 
$(0,1)\times G_{\varepsilon/2}$. Observe that by construction
(see \cite{Kr_12})
$$
 P_{0}(u_{ij})\leq P(u_{ij},u_{i},0)
\leq P(u_{ij},u_{i},u)+K_{1}u_{+}.
$$
It follows that
$$
P_{0}[u_{K}]-K_{1}N-K\leq P[u_{K}]-K
$$
and $u_{K}$ is a solution of \eqref{4.1.2}, which is unique
by the maximum principle. Now, by Theorem 5.4 of
\cite{Kr_14}
$$
[u_{K}]_{C^{1+\chi}(G_{\varepsilon})}\leq
N\varepsilon^{-1-\chi}
$$
and \eqref{3.30.8} for $u_{K}$
 follows from interpolation inequalities.

The function $w=-v_{K}$ is a solution of the equation
$$
\max(-F[-w],P[w]-K)=0,
$$
which is treated similarly to the above.
Observe that in \cite{Kr_14} the operator $F[u]$
is not necessarily given by \eqref{1.16.1}.
The theorem is proved.

The following result is central in this paper.
Fix a constant $\tau\in(0,1)$.

\begin{theorem}
                                            \label{theorem 4.1.2}
For $K\to\infty$ we have $|u_{K}-v_{K}|\to0$
uniformly in $G$. Moreover, if
\begin{equation}
                                                    \label{4.5.1}
\omega(t)=t^{\tau},
\end{equation}
then there exist    constants $\xi\in(0,1)$
depending only on $\tau$,
    $d$, $K_{0}$, and $\delta$
and   $N\in(0,\infty)$,   depending only on  
  $\tau$,  $d$, $K_{0}$,  $\delta$, and   $G$,
such that, if $K\geq N $, then
\begin{equation}
                                                 \label{3.30.01}
|u_{K}-v_{K}|\leq  NK^{-\xi}
\end{equation}
 in $G $.

\end{theorem}

We prove Theorem \ref{theorem 4.1.2}
 in Section \ref{section 4.3.3}

\begin{theorem}
                                            \label{theorem 4.1.3}
(i) The limit
$$
w:=\lim_{K\to\infty}u_{K}
$$
exists,

(ii) The function $w$ is a unique continuous in $\bar G$
viscosity solution
of \eqref{4.1.4} with boundary condition $w=g$ on
$\partial G$,

(iii) If condition \eqref{4.5.1} is satisfied, then
for large enough $K$ we have $|u_{K}-w|\leq  NK^{-\xi}$, 

(iv) For any $\varepsilon\in(0,1]$
(such that $G_{\varepsilon}\ne\emptyset$)
$$
 \|w\|_{C^{1+\chi}(G_{\varepsilon})}\leq
N\varepsilon^{-1-\chi},
$$
where $N$ is the constant from \eqref{3.30.8}.

\end{theorem}

Assertions (i), (iii), and (iv) are simple consequences
of Theorems \ref{theorem 3.31.1}
and \ref{theorem 4.1.2} and the maximum principle.
Indeed,
notice that $F[u_{K}]\leq0$ and $F[v_{K}]\geq0$.
Hence by the maximum principle $u_{K}\geq v_{K}$.
Furthermore, again by the maximum principle
$u_{K}$ decreases and $v_{K}$ increases as $K$
increases.  This takes care
of assertions (i), (iii), and (iv).

Assertion (ii) is proved in Section \ref{section 4.3.4}.

The following result is obtained be referring to the proof
of Theorem 2.1 of \cite{Kr_14_2}
(see Section 5 there) and using 
assertion (iii) of our 
Theorem \ref{theorem 4.1.3}, that was used in
\cite{Kr_14_2} in the case of Lipschitz continuous coefficients.
 To state it we introduce the necessary objects.

As is well known (see, for instance, \cite{KT92}),
there exists a finite set $\Lambda=\{l_{1},...,l_{d_{2}}\}
\subset\bZ^{d}$  containing all vectors from the standard 
orthonormal basis of $\bR^{d}$
such that one has the following representation
$$
L^{\alpha\beta}u(x)=
a^{\alpha\beta}_{k}(x)D_{l_{k}}^{2}u(x)+
\bar b^{\alpha\beta}_{k}(x)D_{l_{k}}u(x)
-c^{\alpha\beta}(x)u(x),
$$
where $D_{l_{k}}u(x)=\langle D u ,l_{k}\rangle$,
$a^{\alpha\beta}_{k}$ and $\bar b^{\alpha\beta}_{k}$ are certain
bounded functions and $a^{\alpha\beta}_{k}\geq\delta_{1}$,
with a constant $\delta_{1}>0$. One can even arrange for
such representation to have the coefficients 
$a^{\alpha\beta}_{k}$ and $\bar b^{\alpha\beta}_{k}$ 
with the same regularity
properties with respect to $x$ as the original ones
$a^{\alpha\beta}_{ij}$ and $b_{i}^{\alpha\beta}$
(see, for instance, Theorem 3.1 in \cite{Kr_11}).
Define
$B$ as the smallest closed ball containing $\Lambda$,
and for $h>0$ set $\bZ^{d}_{h}=h\bZ^{d}$,
$$
G_{(h)}=G\cap\bZ^{d}_{h},\quad
G^{o}_{(h)}=\{x\in \bZ^{d}_{h}:x+hB\in G\},\quad
\partial_{h}G=G_{(h)}\setminus G^{o}_{(h)}.
$$

Next, for $h>0$ we introduce
$$
\delta_{h,l_{k}}u(x)=\frac{u(x+hl_{k})-u(x)}{h} ,
$$
$$
\Delta_{h,l_{k}}u(x)=\frac{u(x+hl_{k})-
2u(x)+u(x-hl_{k})}{h^{2}} ,
$$
$$
L^{\alpha\beta}_{h}u(x)=
a^{\alpha\beta}_{k}(x)\Delta_{h,l_{k}}u(x)+
\bar b^{\alpha\beta}_{k}(x)\delta_{h,l_{k}}u(x)
-c^{\alpha\beta}(x)u(x),
$$
$$
F_{h}[u](x)=
\supinf_{\alpha\in A\,\,\beta\in B}[L^{\alpha\beta}_{h}u(x)
+f^{\alpha\beta}(x)].
$$

It is a simple fact shown, for instance, in \cite{KT92}
that   for each sufficiently small $h$ there exists a unique
  function $w_{h}$ on $G_{(h)}$ such that
$F_{h}[w_{h}]=0$ on $G^{o}_{(h)}$ and 
$w_{h}=0$ on $\partial_{h}G$.

Here is the result we were talking about above.

\begin{theorem}
                                          \label{theorem 4.1.4}
Let condition \eqref{4.5.1} be satisfied
and $g=0$.
Then there exist constants $N$ and $\eta>0$ such that
for all sufficiently small $h>0$ we have on $G_{(h)}$ that
$$
|w_{h}- w|\leq Nh^{\eta}.
$$

\end{theorem}

We imposed the condition $g=0$ in order to be able to
refer directly to the arguments in \cite{Kr_14_2},
where $g=0$. Actually, the result of \cite{Kr_14_2}
can be easily extended to cover the case $g\in C^{1,1}$
(and even go much further),
and then Theorem \ref{theorem 4.1.4} will also cover this case.

\mysection{An auxiliary result}
                                           \label{section 4.3.2}

In the following theorem   
 $G$ can be just any bounded domain.
\begin{theorem}
                                               \label{theorem 3.27.1}
Let   $u,v\in C^{2}(\bar G )$
be such that for a constant $K\geq1$
\begin{equation}
                                                 \label{3.30.5}
\max(F[u],P[u]-K)\geq0\geq\min(F[v],-P[-v]+K)
\end{equation}
in $G $ and $v\geq u$
on $\partial G $. Also assume that, for a constant $M \in[1,\infty)$, 
\begin{equation}
                                                 \label{3.28.06}
\|u,v\|_{C^{1+\chi}(G )}\leq M .
\end{equation}

Then there exist  a  constant    $N \in(0,\infty)$,
depending only on $\tau$, the diameter of $G$,
    $d$, $K_{0}$, and $\delta$,
 and a  constant  $\eta\in(0,1)$,
depending only on $\tau$,
    $d$, $K_{0}$, and $\delta$,
such that, if $K\geq NM^{1/\eta}$, then
\begin{equation}
                                                 \label{3.30.1}
u-v\leq  NK^{-\chi/4}+ NM\omega(M^{-1/\tau}K^{-1})
\end{equation}
 in $G $.

\end{theorem}
\begin{remark}
                                            \label{remark 4.5.2}
Observe that for $\omega=t^{\tau}$ estimate
\eqref{3.30.1} becomes $u-v\leq  NK^{-\chi/4}+K^{-\tau}$.
\end{remark}

As we have mentioned in the Introduction, to prove this theorem,
we are going to adapt to our situation an argument
from Section V.~A of \cite{CIL}.
First we introduce $\psi$ as a global barrier for $G$,
that is $\psi\in C^{2}(G)$, $\psi\geq1$,
$$
a_{ij}D_{ij}\psi+b_{i}D_{i}\psi\leq-1  
$$
in $G$
for any $(a_{ij})\in \bS_{\hat \delta}, |(b_{i})|\leq K_{1}$.
Such a $\psi$ can be found in the form
$  \cosh \mu R-\cosh\mu|  x|$ for sufficiently large $\mu$
and $R$.

Then we take and fix a radially symmetric
  function $\zeta=\zeta( x)$ of class
$ C^{\infty}_{0}(\bR^{d })$ with support in $\{x:|x|<1\}$.
For $\varepsilon>0$ we define $\zeta_{\varepsilon}( x)
=\varepsilon^{-d }\zeta( \varepsilon^{-1}x)$
and for locally summable $u( x)$ introduce
$$
u^{(\varepsilon)}( x)=u( x)*\zeta_{\varepsilon}( x).
$$
Recall some standard properties of these  mollifiers
in which no regularity properties of $G$ are required:
 If $u\in C^{ 1+\chi}(G)$, then with a constant $N$
depending only on $d$  and $\zeta$ (recall that
$\chi$ depends only on $d$ and $\delta$)
$$
\varepsilon^{-1-\chi}|u-u^{(\varepsilon)}|
+\varepsilon^{ -\chi}|Du-Du^{(\varepsilon)}|
\leq N\|u
\|_{C^{ 1+\chi}(G)},
$$
\begin{equation}
                                                    \label{3.31.1}
|u^{(\varepsilon)}
 |+|Du^{(\varepsilon)}| +\varepsilon^{1-\chi}|D^{2}u^{(\varepsilon)}|
 +\varepsilon^{2-\chi}|D^{3}u^{(\varepsilon)}|
\leq N\|u
\|_{C^{ 1+\chi}(G)}
\end{equation}
in 
 $
G_{\varepsilon}$ (introduced in \eqref{4.6.1}).

Define the functions
$$
\bar u=u/\psi,\quad\bar v=v/\psi. 
$$
Increasing $M $ if necessary, we may assume that
\eqref{3.28.06} holds with $\bar u,\bar v$ in place of $u,v$.
This increase, of course, will be affected by the diameter of $G$,
$\delta$,  $d$, and $K_{0}$.

Next, take  constants $\nu,\varepsilon_{0}\in(0,1)$,
recall \eqref{4.14.2},  
introduce 
$$
\varepsilon=\varepsilon_{0} K^{-(1-\gamma)/(2\gamma)},
$$
and
consider the function 
$$
W(x,y)=\bar u(x)-\bar u^{(\varepsilon)}(x)-
[\bar v(y)-\bar u^{(\varepsilon)}(y)]-\nu K|x-y|^{2}
$$
in $  \bar G _{\varepsilon}
\times \bar G _{\varepsilon}$.
Denote by  $(\bar x, \bar y)$ a maximum point
of $W$ in $  \bar G _{\varepsilon}
\times \bar G _{\varepsilon}$.
 Below by $N$  with indices or without them we denote 
various constants
depending only on 
    $d$, $K_{0}$, $\delta$, and  the diameter of $G$,
unless specifically stated otherwise.
By the way recall that $\chi$ depends only on $d$ and $\delta$.

\begin{lemma}
                                              \label{lemma 3.30.1}
There exist  a constant    $\nu\in(0,1)$,
depending only on 
    $d$, $K_{0}$, and $\delta$, and a constant $N$
  such that if
\begin{equation}
                                            \label{3.31.9}
K\geq N\varepsilon_{0}^{(\chi-1)/\eta_{1}}M^{1/\eta_{1}},
\end{equation}
where $  \eta_{1}=1-(1-\chi)(1-\gamma)/(2\gamma)$ ($>0$),
and $ \bar x ,\bar y \in G _{\varepsilon}$, then

(i) We have
\begin{equation}
                                            \label{3.27.6}
2\nu K|\bar x-\bar y|\leq NM \varepsilon^{\chi},
\quad|\bar x-\bar y|\leq\varepsilon/2,
\end{equation}

(ii) For any $\xi,\eta\in\bR^{d}$
\begin{equation}
                                            \label{3.27.7}
D_{ij}[\bar u -\bar u^{(\varepsilon)} ](\bar x)\xi^{i}\xi^{j}
-D_{ij}[\bar v -\bar u^{(\varepsilon)} ](\bar y)\eta^{i}\eta^{j}
\leq 2\nu K|\xi-\eta|^{2},
\end{equation}

(iii) We have
\begin{equation}
                                               \label{3.27.8}
\supinf_{\alpha\in A\,\,\beta\in B}\big[\hat{a}_{ij}^{\alpha\beta}
D_{ij}\bar u+\hat b^{\alpha\beta}_{i}D_{i}\bar u-
\hat c^{\alpha\beta}\bar u+f^{\alpha\beta}\big](\bar x)\geq0,
\end{equation}

\begin{equation}
                                               \label{3.28.1}
\supinf_{\alpha\in A\,\,\beta\in B}\big[\hat{a}_{ij}^{\alpha\beta}
D_{ij}\bar v+\hat b^{\alpha\beta}_{i}D_{i}\bar v-
\hat c^{\alpha\beta}\bar v+f^{\alpha\beta}\big](\bar y)\leq0.
\end{equation}

where
$$
\hat{a}_{ij}^{\alpha\beta}=\psi a_{ij}^{\alpha\beta},
\quad \hat b^{\alpha\beta}_{i}=  b^{\alpha\beta}_{i}\psi
+2a_{ij}^{\alpha\beta}D_{j}\psi,\quad
\hat c^{\alpha\beta}=-L^{\alpha\beta}\psi.
$$

\end{lemma}

Proof. The first inequality in
\eqref{3.27.6}  follows from the fact that the first derivatives
of $W$ with respect to $x$ vanish  at $\bar x$, that is
$D(\bar u-\bar u^{(\varepsilon)}
)(\bar x)=2\nu K(\bar x-\bar y)$.
Also the matrix of second-order derivatives of $W$
is nonpositive at $(\bar x,\bar y)$, which yields (ii).

By taking $\eta=0$ in \eqref{3.27.7} and using the fact that
$|D^{2}\bar u^{(\varepsilon)}|\leq N M
\varepsilon^{\chi-1}$ we see that
$$
D^{2}\bar u(\bar x)\leq 2\nu K+NM
\varepsilon^{\chi-1}.
$$
Furthermore 
$$
D_{ij} u=\psi D_{ij} \bar u+(D_{i}\psi)D_{j}\bar u
+(D_{i}\bar u)D_{j}\psi+(D_{ij}\psi)\bar u,
$$
which implies that
$$
D^{2}  u(\bar x)\leq N(\nu K+ M
\varepsilon^{\chi-1}),\quad P[u](\bar x)
\leq N_{1}(\nu K+ M
\varepsilon^{\chi-1} ).
$$
We now choose and fix $\nu$  such that
\begin{equation}
                                            \label{3.28.2}
N_{1}\nu\leq1/4,\quad   .
\end{equation}
As is easy to see $M\varepsilon^{\chi-1}\leq \nu K$ for $K$
satisfying \eqref{3.31.9} with an appropriate 
$N $.
Then   we have $P[u](\bar x)\leq   K/2<K$, which  
 implies that  
 $F[u](\bar x)\geq0$. We have just proved \eqref{3.27.8}.

Moreover, not only $M\varepsilon^{\chi-1}\leq \nu K$ for
$K$
satisfying \eqref{3.31.9}, but also
$NM\varepsilon^{\chi-1}\leq \nu K$, where $N$ is taken from
\eqref{3.27.6}, if we increase $N$ in \eqref{3.31.9}.
This yields the second inequality in \eqref{3.27.6}.

Similarly, \eqref{3.27.7} with $\xi=0$ implies that
 $
D^{2}\bar v(\bar y)\geq -2\nu K-NM
\varepsilon^{\chi-1}
 $, and, with perhaps different $N_{1}$, that
$P[-\bar v](\bar y)\leq  N_{1}(\nu K+M\varepsilon^{\chi-1})$. 
By denoting by $N_{1}$
the largest of the above $N_{1}$'s, defining $\nu$ by
\eqref{3.28.2}, and taking $K\geq N\varepsilon_{0}^{\chi-1}M^{1/\eta_{1}}$  
 we see that $-P[-\bar v](\bar y) >-K$,
  $F[\bar v](\bar y)\leq0$, and hence
\eqref{3.28.1} holds. The lemma is proved.

{\bf Proof of Theorem \ref{theorem 3.27.1}}.
Fix a (large) constant $\mu>0$ to be specified later
as a constant, depending only on  $d$, $K_{0}$,  $\delta$,
and the diameter of $G$, and
first assume that 
$$
W\leq 2K^{-\chi/4}+\mu M\omega(M^{-1/\tau}K^{-1})
$$
in $  \bar G _{\varepsilon}
\times \bar G _{\varepsilon}$. Observe that for 
any point $x\in G $
one can find a point $y\in G _{\varepsilon}$
with $|x-y|\leq\varepsilon$ and then
$$
\bar u(x)-\bar v(x)\leq \bar u(y)-\bar v(y)
+2M \varepsilon\leq W(y,y)+2M \varepsilon
$$
$$
\leq 2K^{-\chi/4}+\mu M\omega(M^{-1/\tau}K^{-1})
+2M \varepsilon.
$$
It follows that if
\begin{equation}
                                                   \label{3.31.7}
K\geq (NM)^{1/\eta_{2}},
\end{equation}
where $\eta_{2}:=(1-\gamma)/(2\gamma)-\chi/4$ ($>0$), 
then \eqref{3.30.1} holds.

It is clear now that, to prove the theorem, it suffices to find $N$
and $\mu$
such that the inequality  
\begin{equation}
                                                 \label{3.28.5}
W(\bar x,\bar y)\geq 2K^{-\chi/4}
+\mu M\omega(M^{-1/\tau}K^{-1})
\end{equation}
is impossible if $K\geq NM^{1/\eta}$ with $N$
and $\eta$ as in the statement of the theorem.
 Of course, we will argue by 
contradiction and suppose that \eqref{3.28.5} holds.

 Obviously, $W\geq -4M $ in $  \bar G _{\varepsilon}
\times \bar G _{\varepsilon}$  
and at a maximum point $(\bar x, \bar y)$ it holds that
$$
\nu K|\bar x-\bar y|^{2}\leq 8M  .
$$
It follows that (recall that $\nu$ is already fixed)
  $|\bar u^{(\varepsilon)}(\bar x)-
\bar u^{(\varepsilon)}(\bar y)|\leq M |\bar x-\bar y|
\leq N M^{1/2}  K ^{-1/2}$, and we have from \eqref{3.28.5} that
$$
\bar u(\bar x)- 
 \bar v(\bar y) -
\nu K|\bar x-\bar y|^{2}\geq 2 K^{-\chi/4}-N  M^{1/2}  K ^{-1/2}
+\mu M\omega(M^{-1/\tau}K^{-1}),
$$
\begin{equation}
                                            \label{3.27.4}
\bar u(\bar x)- 
 \bar v(\bar y) -
\nu K|\bar x-\bar y|^{2} 
\geq K^{-\chi/4}+\mu M\omega(M^{-1/\tau}K^{-1}),
\end{equation}
where the last inequality holds provided that 
\begin{equation}
                                                   \label{3.31.6}
K\geq N 
M^{1/\eta_{3}}
\end{equation}
 with $\eta_{3}=1-\chi/2>0$.
Here if $\bar x$ or $\bar y$ are on $\partial G _{\varepsilon}$, then
for appropriate $\hat x\in\partial G $ and $\hat y\in\partial G $
either 
$$
\bar u(\bar x)- 
 \bar v(\bar y)\leq M \varepsilon+ \bar v(\hat x)-\bar v(\bar y)
\leq M (2\varepsilon+|\bar x-\bar y|)
$$
or
$$
\bar u(\bar x)- 
 \bar v(\bar y)\leq \bar u(\bar x)-\bar u(\hat y)+
M \varepsilon 
\leq M (2\varepsilon+|\bar x-\bar y|).
$$
In any case
$$
2\varepsilon M +N M^{1/2}(\nu K)^{-1/2}-
\nu K|\bar x-\bar y|^{2} \geq K^{-\chi/4},
$$
  which is impossible for 
\begin{equation}
                                                     \label{3.31.8}
 K\geq N (M^{1/\eta_{2}}
+M^{1/\eta_{3}}).
\end{equation}
Consequently, for such $K$, $(\bar x,\bar y)\in G_{\varepsilon}\times
G_{\varepsilon}$.

Upon combining   \eqref{3.31.7},
\eqref{3.31.6}, \eqref{3.31.8}, and \eqref{3.31.9}
we conclude that there exists $N\in(0,\infty)$
and $\eta_{0}\in(0,1)$
depending only on  
    $d$, $K_{0}$, and $\delta$ such that,
if
\begin{equation}
                                                     \label{3.30.3}
K\geq N\varepsilon_{0}^{(\chi-1)/\eta_{0}}M^{1/\eta_{0}},
\end{equation}
then
\eqref{3.28.5} implies 
\eqref{3.27.4}
and that  $ \bar x ,\bar y \in G _{\varepsilon}$,
 so that
we can use  
the concussions of Lemma \ref{lemma 3.30.1}.

By   denoting $\sigma^{\alpha\beta}=(\hat a^{\alpha\beta})^{1/2}$
we may write
$$
\hat{a}_{ij}^{\alpha\beta}(x)=\sigma^{\alpha\beta}_{ik}
\sigma^{\alpha\beta}_{jk},
$$
and then \eqref{3.27.7} for $\xi^{i}=
\sigma^{\alpha\beta}_{ik}(\bar x)$ and $\eta^{i}
=\sigma^{\alpha\beta}_{ik}(\bar y)$ implies that
$$
\hat{a}_{ij}^{\alpha\beta}(\bar x)D_{ij}\bar u(\bar x)
\leq  \hat{a}_{ij}^{\alpha\beta}(\bar y)D_{ij}\bar v(\bar y)
+\hat{a}_{ij}^{\alpha\beta}(\bar x)D_{ij}
\bar u^{(\varepsilon)}(\bar x)
-\hat{a}_{ij}^{\alpha\beta}(\bar y)D_{ij}
\bar u^{(\varepsilon)}(\bar y)+J,
$$
where
$$
J:=2\nu K\sum_{i,k=1}^{d}|\sigma^{\alpha\beta}_{ik}(\bar x)
-\sigma^{\alpha\beta}_{ik}(\bar y)|^{2}
\leq NK|\bar x-\bar y|^{2\gamma},
$$
and the estimate of $J$ is valid 
because $\hat a^{\alpha\beta}$ is uniformly nondegenerate
and its square root possesses the same smoothness
properties as $\hat a^{\alpha\beta}$. Also note that
$$
\hat{a}_{ij}^{\alpha\beta}(\bar x)D_{ij}\bar u^{(\varepsilon)}(\bar x)
-\hat{a}_{ij}^{\alpha\beta}(\bar y)D_{ij}
\bar u^{(\varepsilon)}(\bar y)\leq
[\hat{a}_{ij}^{\alpha\beta}(\bar x)-
\hat{a}_{ij}^{\alpha\beta}(\bar y)]
D_{ij}\bar u^{(\varepsilon)}(\bar x)
$$
$$
+N
|D^{2}\bar u^{(\varepsilon)}(\bar x)-
D^{2}\bar u^{(\varepsilon)}(\bar y)|
\leq NM|\bar x-\bar y|^{ \gamma}\varepsilon^{\chi-1}
+NM\varepsilon^{\chi-2}|\bar x-\bar y|,
$$
where the last inequality is obtained by the mean-value theorem
relying on the fact that $|\bar x-\bar y|\leq\varepsilon/2$,
so that the straight segment connecting these points
lies inside $G_{\varepsilon/2}$.
Hence, in light of \eqref{3.27.6}
 we get
\begin{equation}
                                                 \label{3.28.4}
\hat{a}_{ij}^{\alpha\beta}(\bar x)D_{ij}\bar u(\bar x)
\leq  \hat{a}_{ij}^{\alpha\beta}(\bar y)D_{ij}\bar v(\bar y)
+NI,
\end{equation}
where
$$
I:=M^{ 2\gamma}K^{1-2\gamma}\varepsilon^{2\gamma \chi}
+M^{1+\gamma}K^{-\gamma}\varepsilon^{\gamma\chi +\chi-1}
+M^{2}K^{-1}\varepsilon^{2\chi-2} =I_{1}+I_{2}+I_{3} .
$$

It turns out that
$$
I_{1}=M^{ 2\gamma}\varepsilon_{0}^{2\gamma\chi}K^{-\chi/4},\quad
I_{2}=M^{1+\gamma}\varepsilon_{0}^{\gamma\chi +\chi-1}K^{-\chi/4-\theta },
$$
$$
I_{3}=M^{2}\varepsilon_{0}^{2\chi-2}K^{-\chi/4-2\theta },
$$
where
$$
\theta=(1-\gamma)(8\gamma)^{-1} \chi >0.
$$
Clearly, 
$$
I\leq  M^{2}\varepsilon_{0}^{2\gamma\chi}K^{-\chi/4}
[1+\varepsilon_{0}^{-\gamma\chi+\chi-1}K^{-\theta}
+\varepsilon_{0}^{2\chi-2-2\gamma\chi}K^{ -2\theta }],
$$
where the expression inside the square brackets is less
than 3 provided that 
\begin{equation}
                                                 \label{3.30.4}
K\geq\varepsilon_{0}^{(\chi-\chi\gamma-1)/\theta},
\end{equation}
which we suppose to hold below.
Then  in light of \eqref{3.28.4} we get  
\begin{equation}
                                                 \label{3.28.6}
\hat{a}_{ij}^{\alpha\beta}(\bar x)D_{ij}\bar u(\bar x)
\leq  \hat{a}_{ij}^{\alpha\beta}(\bar y)D_{ij}\bar v(\bar y)
+N M^{2}\varepsilon_{0}^{2\gamma\chi}K^{-\chi/4}.
\end{equation}

Furthermore,
$$
D_{i}\bar u(\bar x)=2\nu K(\bar x_{i}-\bar y_{i})+
D_{i}\bar u^{(\varepsilon)}(\bar x),\quad
 D_{i}\bar v(\bar y)=
2\nu K(\bar x_{i}-\bar y_{i})
+D_{i}\bar u^{(\varepsilon)}(\bar y),
$$
and the rough estimate (see \eqref{3.27.6})
$$
|\bar x-\bar y|\leq N\varepsilon_{0}^{\chi}MK^{-1}
$$
leads to (cf. Remark \ref{remark 4.5.1}
and \eqref{3.27.6})
$$
\hat b^{\alpha\beta}D_{i}\bar u^{(\varepsilon)}(\bar x)-
\hat b^{\alpha\beta}D_{i}\bar u^{(\varepsilon)}(\bar y)\leq
NM\omega(|\bar x-\bar y|)\leq 
NM\omega(N\varepsilon_{0}^{\chi}MK^{-1}).
$$
Finally,  
$$
f^{\alpha\beta}(\bar x)\leq f^{\alpha\beta}(\bar y)+
N\omega(|\bar x-\bar y|)\leq f^{\alpha\beta}(\bar y)+
N \omega(N\varepsilon_{0}^{\chi}MK^{-1}),
$$
$$
-\bar c^{\alpha\beta}\bar u(\bar x) 
+\bar c^{\alpha\beta}\bar v(\bar y)=\bar c^{\alpha\beta}(\bar x)
[\bar u(\bar x)-\bar v(\bar y)]
$$
$$
+\bar v(\bar y)
[\bar c^{\alpha\beta}(\bar y)-\bar c^{\alpha\beta}(\bar x)]
\leq-[\bar u(\bar x)-\bar v(\bar y)]
+NM\omega(N\varepsilon_{0}^{\chi}MK^{-1}),
$$
where the last inequality follows from the fact that
$\bar c^{\alpha\beta}\geq1$ and $\bar u(\bar x)-\bar v(\bar y)\geq0$
(see \eqref{3.27.4}).

We infer from \eqref{3.27.8}, \eqref{3.28.1},
 and the last estimates that
$$
0\leq \supinf_{\alpha\in A\,\,\beta\in B}\big[ 
\hat{a}_{ij}^{\alpha\beta}D_{ij}\bar v+
\hat b^{\alpha\beta}_{i}D_{i}\bar v-\hat c^{\alpha\beta}
\bar v+f^{\alpha\beta}\big)(\bar y)  \big]
$$
$$
  -[\bar u(\bar x)-\bar v(\bar y)]
+N_{1} M^{2}\varepsilon_{0}^{2\gamma\chi}K^{-\chi/4}
+N_{1}M\omega(N\varepsilon_{0}^{\chi}MK^{-1}),
$$
$$
\bar u(\bar x)-\bar v(\bar y)\leq 
N_{1} M^{2}\varepsilon_{0}^{2\gamma\chi}K^{-\chi/4}
+N_{1}M\omega(N_{1}\varepsilon_{0}^{\chi}MK^{-1}).
$$
We can certainly assume that $N_{1}\geq1$.
Then
we take $\mu=2N_{1}$ and take $\kappa,\xi\in(0,1)$
and $\xi\in(0,\infty)$, depending only on
$\tau$, $K_{0}$, $d$,  $\delta$, and the diameter of $G$, such that
for $\varepsilon_{0}=\xi M^{-1/\kappa}$ and all $M\geq1$
we have 
$$
 N_{1}\varepsilon_{0}^{\chi}M\leq M^{-1/\tau},\quad
N_{1} M^{2}\varepsilon_{0}^{2\gamma\chi}\leq 1/2.
$$
Then 
we arrive at a contradiction with
\eqref{3.27.4} and, 
since now \eqref{3.30.4} and \eqref{3.30.3}
are satisfied if $K\geq NM^{1/\eta}$ for  appropriate
$\eta$ and $N$, the theorem is proved.

\mysection{Proof of Theorem \protect\ref{theorem 4.1.2}}
                                           \label{section 4.3.3}

Fix a sufficiently small $\varepsilon_{0}>0$ such that
$G_{\varepsilon_{0}}\ne\emptyset$ and for $\varepsilon
\in(0,\varepsilon_{0}]$ define
$$
\xi_{\varepsilon,K}:=\max(F[u_{K}^{(\varepsilon)}],
P[u_{K}^{(\varepsilon)}]-K), \quad
\eta_{\varepsilon,K}:=\min(F[v_{K}^{(\varepsilon)}],
-P[-v_{K}^{(\varepsilon)}]+K)
$$
 in $G_{\varepsilon_{0}}$. Since the second-order derivatives
of $u_{K}$ and $v_{K}$ are bounded in $G_{\varepsilon_{0}}$,
we have $ \xi_{\varepsilon,K},\eta_{\varepsilon,K}
\to0$ as $\varepsilon\downarrow0$
in any $\cL_{p}(G_{\varepsilon_{0}})$ for any $K$. Furthermore,
$ \xi_{\varepsilon,K},\eta_{\varepsilon,K}$ are
 continuous. Therefore, there exist smooth functions
$ \bar\xi_{\varepsilon,K},\bar\eta_{\varepsilon,K}$ such that
$$
|\xi_{\varepsilon,K}-\bar \xi_{\varepsilon,K}|+
|\eta_{\varepsilon,K}-\bar\eta_{\varepsilon,K}|\leq\varepsilon,
\quad  \xi_{\varepsilon,K}\leq\bar\xi_{\varepsilon,K},
\quad \eta_{\varepsilon,K}\leq\bar\eta_{\varepsilon,K}
$$
in $G_{\varepsilon_{0}}$.

By Safonov's theorem (see \cite{Sa_84}, \cite{Sa_88}),
for any subdomain $G'$ of $G_{\varepsilon_{0}}$
of class $C^{3}$  there exists a unique
$w_{\varepsilon,K}\in C^{2}(\bar G')$ satisfying
$$
\sup_{a\in\bS_{\hat\delta},|b|\leq K_{1}}[
a_{ij}D_{ij}w_{\varepsilon,K}+b_{i}D_{i}w_{\varepsilon,K}]=-
 (|\bar\xi_{\varepsilon,K}|+|\bar\eta_{\varepsilon,K}|)
$$
in $G'$ with zero boundary condition. Obviously,
$$
\max(F[u_{K}^{(\varepsilon)}-w_{\varepsilon,K}],
P[u_{K}^{(\varepsilon)}-w_{\varepsilon,K}]-K)\geq0,
$$
$$
\min(F[v_{K}^{(\varepsilon)}+w_{\varepsilon,K}],
-P[-v_{K}^{(\varepsilon)}-w_{\varepsilon,K}]+K)\leq0
$$
in $G'$.
After setting
$$
\zeta_{\varepsilon,K}=\sup_{\partial G'}(
u_{K}^{(\varepsilon)}-v_{K}^{(\varepsilon)}-2w_{\varepsilon,K})_{+}
$$
we conclude by Theorem \ref{theorem 3.27.1} applied to
$u_{K}^{(\varepsilon)}-w_{\varepsilon,K}$
and $v_{K}^{(\varepsilon)}+w_{\varepsilon,K}+\zeta_{\varepsilon,K}$
in place of $u$ and $v$, respectively, that
there exist  a  constant    $N \in(0,\infty)$,
depending only on $\tau$, the diameter of $G$,
    $d$, $K_{0}$, and $\delta$,
 and a  constant  $\eta\in(0,1)$,
depending only on $\tau$,
    $d$, $K_{0}$, and $\delta$,
such that, if $K\geq NM_{\varepsilon,K}^{1/\eta}$, then
$$
u_{K}^{(\varepsilon)}-v_{K}^{(\varepsilon)}
\leq  \zeta_{\varepsilon,K}+w_{\varepsilon,K}+
 NK^{-\chi/4}+  M\omega(M^{-1/\tau}K^{-1})
$$
 in $G' $, where $M_{\varepsilon,K}$ is any number satisfying
$$
M_{\varepsilon,K}\geq\|u_{K}^{(\varepsilon)}-w_{\varepsilon,K},
v_{K}^{(\varepsilon)}+w_{\varepsilon,K}+\zeta_{\varepsilon,K}
\|_{C^{1+\chi}(G' )}.
$$

First we discuss
what is happening as $\varepsilon\downarrow0$. By $W^{2}_{p}$-theory
(see, for instance, \cite{Wi09}) 
$w_{\varepsilon,K}\to 0$ in any $W^{2}_{p}$, which by
embedding theorems implies that $w_{\varepsilon,K}\to 0$ in
$C^{1+\chi}(G' )$. Obviously, the constants 
$\zeta_{\varepsilon,K}$ converge in $C^{1+\chi}(G' )$ to
$$
\sup_{\partial G'}(
u_{K} -v_{K}  )_{+}.
$$
Now Theorem \ref{theorem 3.31.1} implies that
for sufficiently small $\varepsilon$ one can take
$N\varepsilon_{0}^{-1-\chi}$ as $M_{\varepsilon,K}$,
where $N$ depends only on $d$, $\delta$, $G$, and $K_{0}$.
Thus for sufficiently small $\varepsilon$,
if $K\geq N\varepsilon_{0}^{-(1+\chi)/\eta}$,
then
$$
u_{K}^{(\varepsilon)}-v_{K}^{(\varepsilon)}
\leq  \zeta_{\varepsilon,K}+w_{\varepsilon,K}+
 NK^{-\chi/4}+N\varepsilon_{0}^{-1-\chi}
\omega(\varepsilon_{0}^{  (1+\chi)/\tau}K^{-1})
$$
in $G'$, which after letting $\varepsilon\downarrow0$
yields
$$
u_{K} -v_{K} 
\leq   
 NK^{-\chi/4}+N\varepsilon_{0}^{-1-\chi}
\omega(\varepsilon_{0}^{ (1+\chi)/\tau}K^{-1})
$$
in $G'$.
The arbitrariness of $G'$ and Lemma \ref{lemma 4.1.1}
now allow us to conclude
that for any   $\varepsilon_{0}>0$, for which
$G_{\varepsilon_{0}}\ne\emptyset$,
$$
u_{K}-v_{K}\leq  NK^{-\chi/4}
+N\varepsilon_{0}^{-1-\chi}
\omega(\varepsilon_{0}^{ (1+\chi)/\tau}K^{-1})
+\sup_{G\setminus G_{\varepsilon_{0}}}(
u_{K} -v_{K}  )_{+}
$$
\begin{equation}
                                                  \label{4.6.2}
\leq  NK^{-\chi/4}
+N\varepsilon_{0}^{-1-\chi}
\omega(\varepsilon_{0}^{(1+\chi)/\tau}K^{-1})+N\varepsilon_{0}
\end{equation}
in $G$ provided that  
\begin{equation}
                                                  \label{4.6.3}
K\geq N\varepsilon_{0}^{-(1+\chi)/\eta}.
\end{equation}
This obviously proves the 
first assertion of the theorem because as is noted in the proof
of Theorem \ref{theorem 4.1.3} we have $v_{K}\leq u_{K}$.

To prove the second assertion observe that for $\omega=
t^{\tau}$ and $\varepsilon_{0}=K^{-\eta/(1+\chi)}$
condition \eqref{4.6.3} becomes $K\geq N$ and \eqref{4.6.2}
becomes
$$
u_{K}-v_{K}\leq  N K^{-\chi/4}+N K^{-\tau} +N
K^{-\eta/(1+\chi)}.
$$
This yields  the desired result and proves the theorem.

\mysection{Proof of assertion (ii)
of Theorem \protect\ref{theorem 4.1.3}}
                                           \label{section 4.3.4}

First we show uniqueness. Let $v$
be a continuous in $\bar G$
 viscosity solution of $F[v]=0$ with boundary data
$g$. Observe that in the notation
from Section \ref{section 4.3.3} we have
$$
F[u_{K}^{(\varepsilon)}+w_{\varepsilon,K}+\kappa \psi]<0
$$
in $G'$ for any $\kappa>0$.
 This and the definition of viscosity solution
imply that the minimum of
$u_{K}^{(\varepsilon)}+w_{\varepsilon,K}+\kappa \psi
-v$ in $\bar G'$ is either positive or is attained
on $\partial G'$. The same conclusion holds after letting
$\varepsilon,\kappa\downarrow0$ and replacing $G'$
with $G_{\varepsilon_{0}}$. Hence, in $G$
$$
u_{K}-v\geq -\sup_{G\setminus G_{\varepsilon_{0}}}
|u_{K}-v|,
$$
which after letting $\varepsilon_{0}\downarrow0$
and then $K\to\infty$
yields $v\leq w$. By comparing $v$ with $v_{K}$ we get
$v\geq w$, and hence uniqueness.

To prove that $w$ is a viscosity solution we need
a lemma, which is an elliptic analog of
Lemma 6.1 of \cite{Kr_14}.  Introduce
$$
F_{0}(u_{ij},x)=F(u_{ij},Dw(t,x),w(t,x),x).
$$
\begin{lemma}
                                           \label{lemma 9.20.1}  
There is a constant $N$ depending only on $d$ and
$\delta$ such that
for any ball $B_{r}$ of radius $r$  with closure in $G$  and
$\phi\in W^{ 2}_{d }(B_{r})
\cap C(\bar B_{r})$ we have on $B_{r}$ that
\begin{equation}
                                             \label{9.20.1}
w\leq \phi+Nr \|( 
F_{0}[\phi])^{+}\|_{L_{d }(B_{r})}
+\max_{\partial B_{r}}(w-\phi)^{+} .
\end{equation}
\begin{equation}
                                             \label{9.20.2}
w\geq \phi-Nr \|( 
F_{0}[\phi] )^{-}\|_{L_{d }(B_{r})}
-\max_{\partial B_{r}}(w-\phi)^{-} .
\end{equation}
\end{lemma}

Proof. Observe that
$$
 -\max(F_{0}[\phi] ,P[\phi]-K )=\max(F_{0}[u_{K}] ,P[u_{K}]-K)
 -\max(F_{0}[\phi] ,P[\phi]-K )
$$

$$
+I_{K}=
  a_{ij}D_{ij}(u_{K}-\phi)
+I_{K},
$$
where $a=(a_{ij})$ is an $\bS_{\hat{\delta}}$-valued
function and
$$
I_{K}=\max(F [u_{K}] ,P[u_{K}]-K )-
\max(F_{0}[u_{K}] ,P[u_{K}]-K ).
$$

It follows by Theorem 5.2  of \cite{Kr78} or Theorem 3.3.11
of \cite{Kr87}
that  
$$
u_{K}\leq \phi+\max_{\partial B_{r}}(u_{K}-\phi)^{+}
$$
\begin{equation}
                                             \label{9.20.3}
+Nr \|( I_{K}+
\max(F_{0}[\phi] ,P[\phi]-K ))^{+}\|_{L_{d }(B_{r})},
\end{equation}
where the constant $N=N(d,\delta)$. Actually the above references
only say that \eqref{9.20.3} holds with $N=N(r,d,\delta)$
in place of $Nr $.
However, the way this constant depends on $r$ is easily discovered
by using   dilations.

Notice that  $u_{K}\to w$ and $Du_{k}\to Dw$ as $K\to\infty$
uniformly in $B_{r}$. Hence $I_{K}\to0$ as $K\to\infty$
uniformly in $B_{r}$.

After that we obtain \eqref{9.20.1} from \eqref{9.20.3}
by letting $K\to\infty$.
In the same way \eqref{9.20.2} is established by considering $v_{K}$.
The lemma is proved.

Now let $\phi\in C^{2}(G)$ and suppose that 
$w-\phi$ attains a local
maximum at $ x_{0}\in G$. Without losing generality
we may assume that $x_{0}=0$, $w(0)-\phi(0)=0$.
Then for all small $r>0$ and balls $B_{r}$ centered at
$x_{0}$ and $\varepsilon>0$ by applying 
\eqref{9.20.1} at the origin
 to $\phi(x)-\varepsilon(r^{2}-|x|^{2})$
in place of $\phi$ we get
$$
\varepsilon r^{2}\leq 
Nr \|( 
F_{0}[\phi -\varepsilon(r^{2}-|\cdot|^{2})])^{+}\|_{L_{d }(B_{r})}.
$$
It follows that
$$
\sup_{B_{r}}[F_{0}[\phi -\varepsilon(r^{2}-|\cdot|^{2})]\geq0,
$$
and by letting first $r\downarrow0$ and then $\varepsilon
\downarrow0$ we conclude that 
$$
0\leq F_{0}[\phi](0)=F(D_{ij}\phi(0),D_{i}\phi(0),\phi(0),0),
$$
where the equality follows from the fact that at $0$
the derivatives of $w-\phi$ vanish.
We have  just proved that $w$ is a viscosity subsolution.

Similarly by using \eqref{9.20.2} one
proves that $w$ is a viscosity supersolution.
This proves the theorem.

\end{document}